# Economic Sizing of Distributed Energy Resources for Reliable Community Microgrids

Chen Yuan and Mahesh S. Illindala
Department of Electrical and Computer Engineering
The Ohio State University
Columbus, Ohio, USA
yuan.261@osu.edu

*Abstract*—Community microgrids offer many advantages for power distribution systems. When there is an extreme event happening, distribution systems can be seamlessly partitioned into several community microgrids for uninterrupted supply to the end-users. In order to guarantee the system reliability, distributed energy resources (DERs) should be sized for ensuring generation adequacy to cover unexpected events. This paper presents a comprehensive methodology for DERs selection in community microgrids, and an economic approach to meet the system reliability requirements. Algorithms of discrete time Fourier transform (DTFT) and particle swarm optimization (PSO) are employed to find the optimal solution. Uncertainties of load demand and renewable generation are taken into consideration. As part of the case study, a sensitivity analysis is carried out to show the renewable generation impact on DERs' capacity planning.

*Index Terms*--Capacity planning, distributed power generation, microgrids, optimization, power systems reliability.

## I. INTRODUCTION

Nowadays, power outages caused by extreme weather and overloading are happening more frequently. To ensure the maximum extent of uninterruptible power supply for end-users, distribution systems could be seamlessly partitioned into community microgrids and each microgrid should be able to operate independently. But one of the challenges is how to guarantee the reliability of the microgrid islanded operation. Most of works about the capacity planning of DERs have been focusing on cost minimization [1], [2]. But very few considered the failure rates of generators and stochastic natures of load and renewable resources towards the system reliability. A study of microgrid generation adequacy indicates that, with a certain generation capacity, more distributed energy resources (DERs) could lead to higher system reliability [3]. The reliability centered generation capacity planning was conducted in [4], but the intermittency and uncertainty of renewable generation were neglected. In this paper, selection of DERs for community microgrids based on both quantitative and qualitative evaluation is presented. Then, an optimal and economic DERs sizing scheme for reliable community microgrids is formulated, with consideration to load uncertainty and renewables unpredictability.

## II. DISTRIBUTED ENERGY RESOURCES SELECTION AND COMMUNITY MICROGRID

### A. Distributed Energy Resources Selection

*1) Levelized cost of energy − Quantitative Assessment:*
Levelized cost of energy (LCOE) is one of the utility's principal metrics for measuring the cost of electricity produced by a generator. It calculates the value of the unit's annualized total cost divided by its estimated annual energy output, as expressed in (1). The parentheses indicate that when the generator is a fuel-powered resource, the fuel cost is accounted in (1); but if it is a renewable resource, a tax credit payback is accounted in (1). The annualized total cost consists of annualized capital cost, $C^{capital}$, operation and maintenance (O&M) cost, $C^{O\&M}$, fuel cost, $C^{fuel}$, and also renewable energy tax incentive payback, $C^{TaxCredit}$, as shown in (2) – (5). The O&M cost includes a fixed part, $C^{fixed\ O\&M}$, and a variable portion, $C^{variable\ O\&M}$. The fixed O&M cost depends on the power rating of each DER, while the variable one is determined by the real operation and related to the energy output. In this paper, the tax incentive is the production tax credit (PTC) — a federal tax incentive that provides financial support for the development of renewable energy [5].

$$LCOE = \frac{C^{capital} + C^{O\&M} + (C^{fuel}) - (C^{TaxCredit})}{P_R \cdot 8760\ h/yr \cdot CF} \quad (1)$$

$$C^{capital} = Overnight\ Capital\ Cost \cdot \frac{r \cdot (1+r)^y}{(1+r)^y - 1} \quad (2)$$

$$C^{O\&M} = C^{fixed\ O\&M} + C^{variable\ O\&M} \quad (3)$$

$$C^{fuel} = F \cdot \sum_{fuel}\left(\frac{\sum_i P[i] \cdot T}{\eta[i]}\right) \cdot Heat\ Rate \cdot LF \quad (4)$$

$$C^{TaxCredit} = PTC \cdot \sum_{re}\left(\sum_i P_{re}[i] \cdot T\right) \quad (5)$$

where, *r* indicates the discount rate referring to the interest rate used in discounted cash flow analysis for setting the present value of future cash flows, *y* is the number of years in a lifetime, *F* indicates the fuel price, $P[i]$ is the power output at time *i*, $\eta[i]$ is the energy conversion efficiency at time *i*, *T* is the sampling time, *LF* is the levelizing factor to estimate future change of fuel cost, and $P_{re}[i]$ is the renewable resources power output at time *i*.

As a financial tool, LCOE is very valuable for the comparison of various generation units [6]. A low LCOE indicates a low cost of electricity generation. For a conventional power plant, the future fuel price is uncertain and largely depending on external factors, while a renewable energy resource has zero fuel cost, although the initial capital cost is high. Besides, governments have policies to encourage the integration of renewable energy resources, like subsidies, tax incentives, feed-in tariff, net-metering program, renewable portfolio standards, etc. Figure 1 is the comparison of LCOE vs. capacity factor among different DERs.

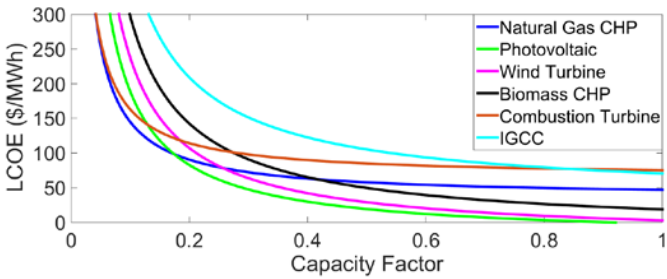

Figure 1. Curves of LCOE vs. capacity factor for DERs

*2) Qualitative Function Deployment — Qualitative Evaluation:* After the LCOE calculation, a qualitative evaluation could further assess DERs' soft indices [6]. This can be made using qualitative function deployment (QFD). As seen in Table I, biomass gensets, natural gas gensets, PV panels, and wind turbines are more suitable options than other choices. Besides, the natural gas is very efficient and has ample supplies in the U.S.

*B. Community Microgrids*

Based on the quantitative comparison and qualitative evaluation, the suitable DERs for community microgrids have been selected out. This subsection displays the community microgrid development within an existing distribution system, as shown in Figure 2. It is formed by integrating local DERs, which are strategically placed near critical loads [6], [7].

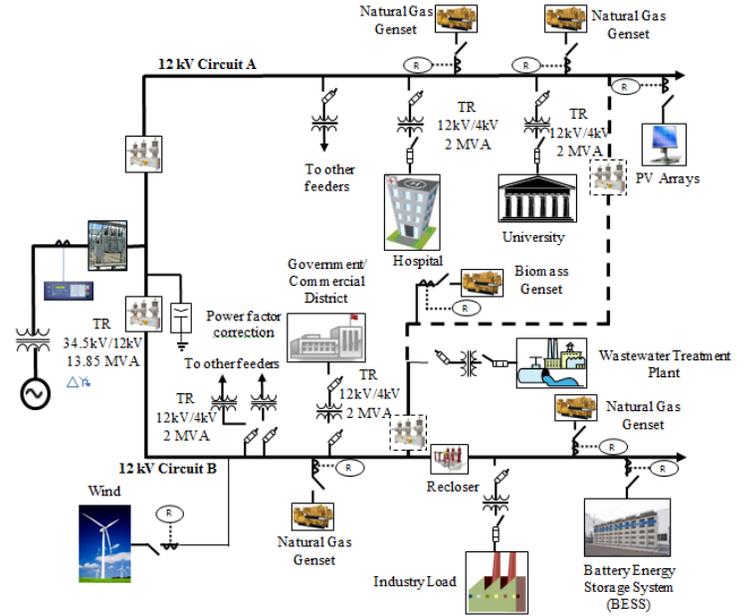

Figure 2. Single-line diagram of a community microgrid

### III. PLANNING RESERVE MARGIN

Unexpected net load changes could cause system instability, overloading, and even power outages. In order to ensure a reliable community microgrid, it is necessary to have adequate reserve margin.

*A. Planning Reserve Margin*

Since it is very difficult to accurately forecast future electricity usage, power losses and generation from renewable energy resources, the concept of planning reserve margin (PRM) is employed to maintain systems reliability. It is a key metric that measures the flexibility to meet customer demands and the ability to handle the loss of system components. The calculation of PRM is coupled with probabilistic analysis to identify the resource adequacy and find out whether the planning capacity is large enough to cover peak load demand, loss of one generation unit, and also uncertainties from load and renewables. It is expressed below.

$$\begin{aligned}Planning\ Capacity &= Peak\ Load \cdot (1 + PRM) \\ &= Peak\ Load + Largest\ Dispatchable\ Generation \\ &\quad Unit + Uncertainties\end{aligned} \quad (6)$$

*B. Impact of Planning Reserve Margin on Reliability Metrics*

System reliability and resource adequacy are not readily observable. For example, one cannot quickly evaluate a system's reliability, like loss of load expectation (LOLE) and

| DERs Options / Customer Requirements | Importance (1-5) | PV Panel | Wind Turbine | Biomass Genset | Natural Gas Genset | Natural Gas Combustion Turbine | Coal–Fired Power Plant (Base-line) |
|---|---|---|---|---|---|---|---|
| LCOE | 5 | 9 | 9 | 3 | 3 | 1 | 1 |
| $CO_2$ Emission Reduction | 5 | 3 | 3 | 9 | 9 | 1 | 0 |
| Fuel Consumption Savings | 4 | 9 | 9 | 9 | 3 | 1 | 0 |
| Outage Time Reduction | 5 | -3 | -3 | 1 | 3 | 3 | 3 |
| Dispatchability | 4 | -1 | -1 | 1 | 3 | 3 | 1 |
| Equipment Lifetime | 3 | 3 | 3 | 1 | 1 | 1 | 3 |
| Comply with the U.S. DOE Target | 5 | 9 | 9 | 9 | 1 | 1 | 0 |
| **Absolute Target** | | **131** | **131** | **153** | **107** | **49** | **33** |

system average interruption frequency index (SAIFI), by simply taking a look at the reserve margin. Based on the probabilistic analysis and Monte Carlo simulation, the general relation between PRM and reliability metrics are plotted in Figure 3 and Figure 4. Situations like unexpected load changes, renewable generation fluctuation, and loss of one generation are covered. Power outages caused by external conditions like bad weather, grounding fault, and distribution line disconnection are out of scope since they can be hardly improved by PRM. Figure 3 displays the LOLE vs. PRM curves with different proportions of the largest dispatchable generator (PLG). Besides, SAIFI vs. PRM curves are presented in Figure 4. Both figures reflect that, when the largest dispatchable generation unit takes up a larger portion of total generation capacity, the system reliability is worse, needing a larger PRM to achieve the same level of reliability.

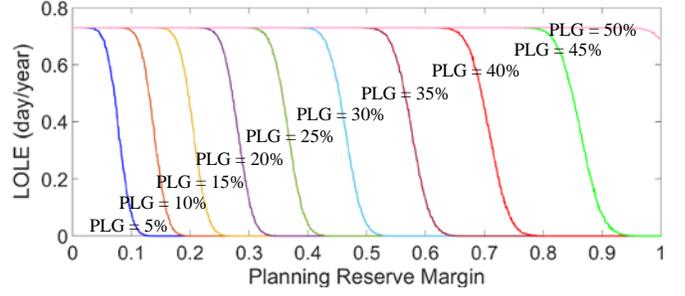

Figure 3. Curves of LOLE vs. planning reserve margin with different proportions of the largest dispatchable generator

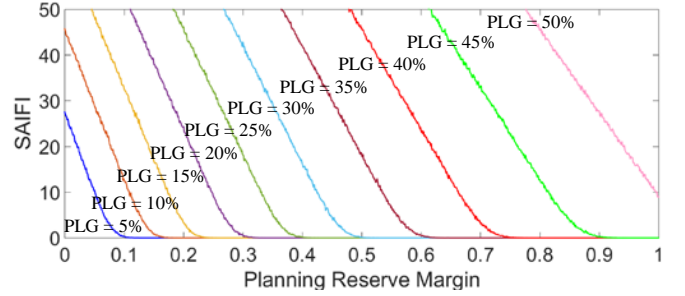

Figure 4. Curves of SAIFI vs. planning reserve margin with different proportions of the largest dispatchable generator

## IV. CAPACITY PLANNING FOR DISPATCHABLE UNITS

### A. Problem Formulation

Renewable energy resources are considered as negative loads and sized to meet the customer needs. So the power capacities of renewables like PV and WT, are determined first per customers requirements. The sizing of gensets and BESS will be processed then. However, with high penetration of renewables, there are some extra challenges because of their uncertainty. For example, a key question is whether the capacity planning for gensets and BESS should take the forecasted generation from renewables into consideration or not? If yes, is it better to accommodate a portion of the renewable generation? These questions will be explored further in this section.

As described in the previous section, a larger reserve margin will lead to a more reliable system. But it will also result in lower efficiency and higher cost. This is because when the total generation capacity is larger, the operation efficiency could be lower with more reserve margin, resulting in higher capital cost, O&M cost, and fuel cost. So, there is a tradeoff between cost and reliability.

Based on the discussion, system reliability is the primary goal of the sizing problem. However, the cost cannot be ignored while carrying out the gensets and BESS sizing. Therefore, the total cost minimization is set as the objective and system reliability requirement as a constraint. By this way, the reliability requirements are satisfied before achieving the cost minimization. So the sizing problem is formulated as follows.

$$\text{Minimize } TC \quad (7)$$

subject to

$$P_{available\ supply}(t) \geq P_l(t) \quad (8)$$

$$LOLE \leq LOLE_{thr} \quad (9)$$

$$PLG = \frac{largest\ dispatchable\ generator}{P_{total}} \leq PRM \quad (10)$$

where,

$$TC = C_{PV} + C_{WT} + C_{Biomass} + C_{NG} + C_{BESS} \quad (11)$$

$$C_{XX} = C^{capital} + C^{O\&M} + (C^{fuel}) - (C^{TaxCredit}) \quad (12)$$

where, $TC$ is the total annualized cost, including PV, WT, biomass genset, natural gas genset and BESS, $P_{available\ supply}(t)$ is the available power supply at time $t$, $P_l(t)$ is the load demand at time $t$, $LOLE_{thr}$ is the threshold value, $P_{total}$ is the total power capacity of DERs, and $C_{XX}$ denotes the annualized cost including capital cost, O&M cost, fuel cost (if the energy resource is fuel-powered), and tax credit payback (if it is a renewable resource). The subscript *XX* points to the type of energy resource — PV, WT, biomass genset, natural gas genset, or BESS.

As indicated in (7) − (10), the problem of sizing gensets and BESS is to minimize the total cost with the satisfaction of load demand and system reliability. In (8), it expresses the necessity to have adequate power supply. The system reliability requirement is presented in (9). When system LOLE is smaller than this value, acceptable system reliability is guaranteed. Furthermore, based on the previous discussion in Section III.B, LOLE will also be affected by the PRM and the PLG of the total available power supply. The PLG should be less than PRM, as shown in (10). Otherwise, the system reliability cannot be satisfied, since the customer interruption is a sure event when the largest dispatchable generator is offline. This is also well supported by Figure 3 and Figure 4. The total annualized cost in (11) covers all DERs' costs. For each type of DER, the annualized cost consists of the annualized capital cost, O&M cost, fuel cost and the tax credit in a year, as presented in (12). Brackets indicate that fuel cost only applies to fossil fuel based generation units and only renewable energy generators have tax credits as payback. Each DER's annualized capital cost, O&M cost, fuel cost and tax credit payback can be calculated using (2) − (5).

## B. Optimization Algorithm

Based on the formulated problem in the Section IV.A, gensets are sized together with BESS to share the net load, $P_{nl}$, and provide adequate reserve margin. The net load can be divided into two parts: a) component with large power that varies smoothly over longer duration, and b) small but frequently fluctuating power component. Gensets could take the smooth (i.e., flat) power variation and the BESS can compensate the small and frequent changes, as shown in (13). $P_{BESS\_AC}(t)$ indicates BESS power output after the dc/ac converter, meaning the power on the ac side of the converter.

$$P_{nl}(t) = \Sigma P_{genset}(t) + P_{BESS\_AC}(t) \qquad (13)$$

Methods of discrete time Fourier transform (DTFT) and particle swarm optimization (PSO), which is a population based stochastic optimization technique, are employed to find the optimal power allocation between gensets and BESS to minimize annualized cost with the satisfaction of the system reliability requirement. As described in Section IV.A, loss of one generation should also be covered within planning reserve margin to obey "N-1" criterion.

The sizing scheme is illustrated in the flowchart shown in Figure 5. As seen in this figure, with the input of power capacities of PV and WT, the net load profile can be achieved based on stochastic models of load and renewable units. The DTFT is applied to obtain the frequency spectrum of the net load. The frequency range of the spectrum is determined by the sampling rate of the net load profile. Therefore, the net load profile in the time domain is converted into components in the frequency domain. Then a randomly initialized cut-off frequency divides the net load into two parts. The low-frequency part of the net load is assigned to gensets, while BESS takes care of the high-frequency power components. This helps lower system's capital cost, since the BESS, in terms of power and energy, will not be oversized. Once the power share for gensets is achieved, based on the initial cut-off frequency, the inverse DTFT is employed to get the power share for gensets in the time domain. However, the process of inverse DTFT may produce negative values. But the gensets cannot absorb power. So the power supply from genset should be updated, as presented in (14). $\Sigma P_{genset}^{update}[k]$ is the updated total power output of gensets at sample point $k$.

$$\Sigma P_{genset}^{update}[k] = \begin{cases} \Sigma P_{genset}[k], & \Sigma P_{genset}[k] \geq 0 \\ 0, & \Sigma P_{genset}[k] < 0 \end{cases} \qquad (14)$$

After making such an adjustment, the BESS power share can be determined by subtracting the power allocation of gensets from the net load. The next step is to size gensets and BESS. As shown in (15), the power capacity of gensets needs to meet the maximum power output and also include a reserve margin to withstand forecast errors and unexpected events. BESS is supposed to have the same conversion efficiency, $\eta_{conv}$, in both charging and discharging. So the BESS power capacity is calculated in (16) and (17). It is discharging when $P_{BESS\_AC}[k] \geq 0$, and charging when $P_{BESS\_AC}[k] < 0$. So the BESS real dc power output before the converter is updated by (16) with consideration of power conversion efficiency. Then its nominal power rating could be determined with (17). Furthermore, the BESS capacity in energy is determined in (18) and (19). The change in stored energy from the original status to $k$th sample point is given by (18). Equation (19) presents the sizing for BESS energy capacity, where $SOC_{max}$ and $SOC_{min}$ are predetermined maximum and minimum values of state of charge (SOC).

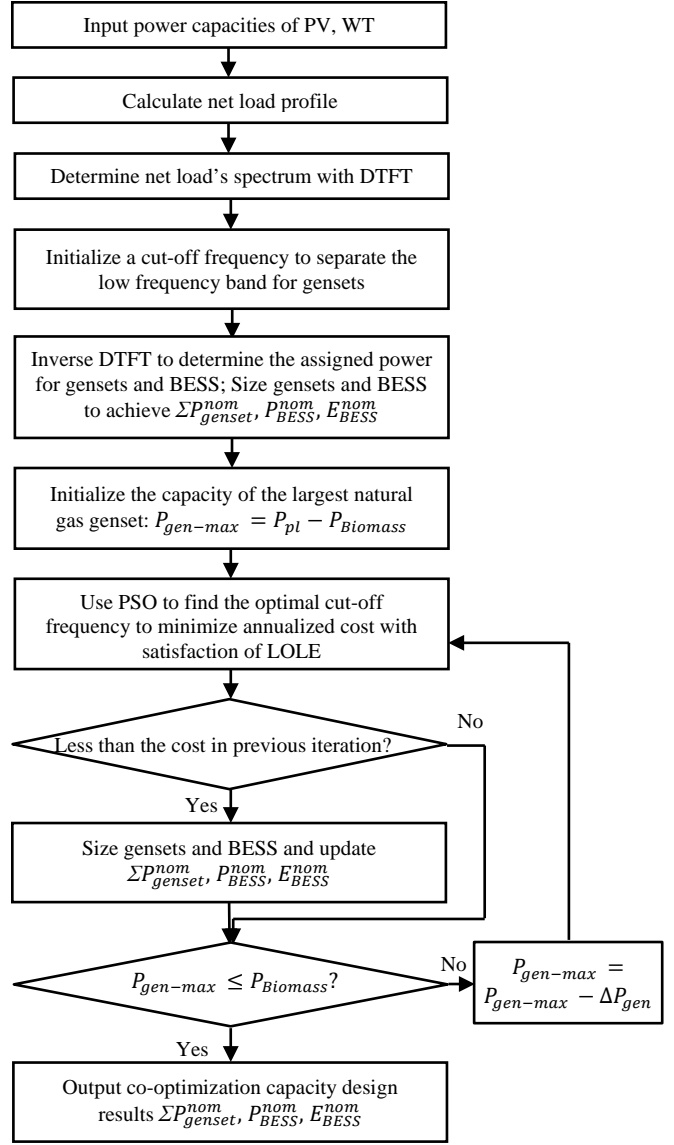

Figure 5. Flowchart of gensets and BESS capacity planning

$$\Sigma P_{genset}^{nom} = max\{\Sigma P_{genset}^{update}[k]\} \cdot (1 + PRM) \qquad (15)$$

$$P_{BESS\_DC}[k] = \begin{cases} P_{BESS\_AC}[k]/\eta_{conv}, & P_{BESS\_AC}[k] \geq 0 \\ P_{BESS\_AC}[k] \cdot \eta_{conv}, & P_{BESS\_AC}[k] < 0 \end{cases} \qquad (16)$$

$$P_{BESS}^{nom} = max\{|P_{BESS\_DC}[k]|\} \qquad (17)$$

$$E[i] = \sum_{k=0}^{i}(P_{BESS\_DC}[k] \cdot T), \qquad i = 0, \dots, N_S \qquad (18)$$

$$E_{BESS}^{nom} = \frac{max\{E[i]\} - min\{E[i]\}}{SOC_{max} - SOC_{min}} \qquad (19)$$

where, $\Sigma P_{genset}^{nom}$ is the nominal power of genset, *PRM* is the portion of reserve margin, $P_{BESS\_DC}[k]$ is BESS power output on the dc side of the power converter at sample point $k$, $P_{BESS}^{nom}$ is the nominal power of BESS, $E[i]$ is the energy difference of BESS from the beginning to sample point $i$, and $E_{BESS}^{nom}$ is the nominal energy capacity of BESS.

After the preliminary sizing, the power capacities of gensets and BESS are achieved. Power share for natural gas gensets equals the subtraction between total power capacity in gensets and biomass genset power capacity. However, since the cut-off frequency was initialized randomly, this may not guarantee the optimal power capacity planning. Therefore, the PSO is used for determining the cut-off frequency to achieve the minimum annualized cost. PSO begins with initialized random solutions and searches for optima by updating iterations within the problem space.

In addition, the largest dispatchable generator has an impact on system's reliability and PRM determination. Therefore, before the PSO, the power capacity of the largest natural gas genset has to be determined. In Figure 5, the power capacity of the largest natural gas genset, $P_{gen-max}$, is initialized as the subtraction of the peak load, $P_{pl}$, and the biomass genset capacity, $P_{Biomass}$. In each iteration, the largest natural gas genset's power capacity is reduced by $\Delta P_{gen}$ until it is smaller than the biomass one. Besides, PSO is implemented in every round to achieve the optimal solution. After all iterations, the optimal solution will be found.

## V. CASE STUDY AND SENSITIVITY ANALYSIS

The proposed strategy for sizing gensets and BESS is implemented and analyzed for the community microgrid in Figure 2, which has 4 MW peak load, 3 MW PV system, 1 MW wind turbine, and 0.5 MW biomass genset. The load data is generated from the load stochastic model based on two years' historical information provided by the local utility, AEP Ohio. The PV and wind output are estimated from their own stochastic models, which are developed using two years' data from references [8], [9]. In addition, $LOLE_{thr}$ is set as 1 day in 10 years. The life time, $y$, is assumed as 20 years and the discount rate, $r$, is 5%.

Table I and Figure 6 present the comparison of minimum annualized cost and the minimum total capacity of gensets for six scenarios. It can be observed that the annualized cost and the gensets total capacity are minimum for 80% renewable generation. Besides, the scenario with 90% renewable generation also has lower minimum annualized cost and gensets total capacity than that for 100% renewable generation.

TABLE I. COMPARISON OF DIFFERENT SCENARIOS

| Scenarios | Annualized Cost ($/MW-year) | Total Capacity of Gensets (MW) | BESS Capacity Power (MW) | BESS Capacity Energy (MWh) | Largest Natural Gas Genset (MW) |
|---|---|---|---|---|---|
| No Renewables | 476,440 | 4.6799 | 0.8507 | 1.1839 | 0.50 |
| 20% Renewables | 460,540 | 4.2330 | 0.8507 | 1.1581 | 0.5 |
| 50% Renewables | 439,840 | 3.6922 | 0.8507 | 0.9908 | 0.5 |
| **80% Renewables** | **427,450** | **3.3019** | **0.8693** | **1.0410** | **0.5** |
| 90% Renewables | 427,580 | 3.3058 | 0.8693 | 1.0410 | 0.5 |
| 100% Renewables | 430,720 | 3.3457 | 0.8955 | 1.0308 | 0.5 |

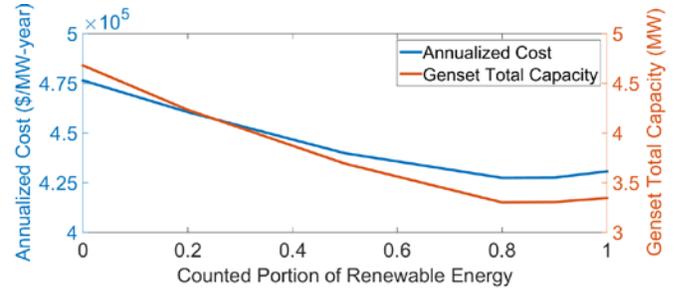

Figure 6. Tendencies of minimum annualized cost and minimum genset total capacity with counted portion of renewable energy

## VI. CONCLUSION

This paper presented an economic sizing method of gensets and BESS for reliable community microgrids. At first, the LCOE based quantitative assessment and QFD based qualitative evaluation were undertaken for various types of DERs to select suitable ones for community microgrids. The sizing scheme was further elaborated for gensets and BESS to ensure system reliability under uncertainties. The employed optimization methodology is based on DTFT and PSO. In the case study, a sensitivity analysis has been conducted to demonstrate that a small margin of renewable generation could provide cost savings and lower capacity of gensets.